\documentclass[draft]{article}
\usepackage{amssymb,amsmath}
\usepackage{graphics}

\input xy
\xyoption{all}
\begin{document}
\title{On the ellipticity of symplectic twistor complexes}
\author{Svatopluk Kr\'ysl \footnote{{\it E-mail address}: krysl@karlin.mff.cuni.cz}\\ {\it \small  Charles University of Prague,  Praha 8,  Czech Republic}
\thanks{The author of this article was supported by the grant
GA\v{C}R  201/08/0397 of the Grant Agency of Czech Republic. The work is a part of the 
research project MSM 0021620839 financed by M\v{S}MT \v{C}R.}}

\maketitle \noindent
\centerline{\large\bf Abstract}  
 
 For a Fedosov manifold (symplectic manifold equipped with a symplectic torsion-free
affine connection $\nabla$) admitting a metaplectic structure, we shall 
investigate two sequences of first order differential operators acting on 
sections of certain bundles over this manifold. The operators are symplectic analogues of the twistor operators known from Riemannian spin 
geometry. Therefore we call the mentioned sequences symplectic twistor sequences.
These sequences are complexes if the connection $\nabla$ is of  Ricci type. 
We shall prove that the so called truncated parts of these complexes are elliptic. This establishes a background for a future analytic study.
 
{\it Math. Subj. Class.:} 22E46, 53C07, 53C80, 58J05

{\it Key words:} Fedosov manifolds, Segal-Shale-Weil representation, Kostant's spinors, elliptic complexes

\section{Introduction}	

  In this article, we prove the ellipticity of certain parts of the so called symplectic twistor complexes. These parts will be defined later in this text. 
The symplectic twistor complexes are two sequences of first order differential operators  defined over Ricci type Fedosov manifolds admitting a metaplectic structure. From reasons clarified bellow in this paper,  we shall call these complexes  left and right symplectic twistor complexes. The mentioned parts of these complexes will be called left and right truncated symplectic twistor complexes, respectively.

 
 
 Now, let us say a few words about the Fedosov manifolds. 
Formally speaking, a Fedosov manifold is a triple $(M^{2l},\omega, \nabla)$ where	
$(M^{2l},\omega)$ is a $2l$ dimensional symplectic manifold and $\nabla$ is an affine torsion-free  symplectic connection.
 By torsion-free and symplectic, one 
means $T^{\nabla}(X,Y):=\nabla_X Y -\nabla_Y X - [X,Y] = 0$ for all vector fields $X,Y \in \mathfrak{X}(M),$ and $\nabla \omega = 0,$ respectively. Connections satisfying these two properties are usually called {\it Fedosov connections} in the honor of Boris Fedosov who used them to obtain a  deformation quantization for symplectic
 manifolds. (See Fedosov \cite{Fedosov}.)  Let us also mention that in  
contrary to torsion-free Riemannian connections, the Fedosov connections are not unique. We refer an interested reader to Tondeur \cite{Tondeur} and Gelfand, Retakh, Shubin \cite{GRS} for more information.

  To formulate the result on the ellipticity of the truncated symplectic twistor complexes, one should know some basic facts on the structure of the curvature of a Fedosov connection. 
In Vaisman \cite{Vaisman}, one can find a proof of a theorem saying that the curvature tensor field of a Fedosov connection splits into 
two parts if $l\geq 2,$ namely into the symplectic Ricci and symplectic 
Weyl curvature tensor fields.  If $l=1,$ only the symplectic Ricci curvature tensor field occurs. 
 Fedosov manifolds with zero symplectic Weyl curvature tensor fields are usually called of Ricci type. (See Vaisman \cite{Vaisman} and Cahen, Schwachh\"ofer \cite{CahenSchw}.)

 After introducing the underlying geometric structure, let us start describing the fields on which 
the differential operators from the symplectic twistor complexes act. The fields are certain exterior  
differential forms with values in the so called symplectic spinor bundle.  The symplectic spinor bundle is an associated
 vector bundle  to the metaplectic bundle. We shall introduce the {\it metaplectic bundle} briefly now.
Because the first homotopy group 
of the symplectic 
group $Sp(2l,\mathbb{R})$ is isomorphic to $\mathbb{Z},$ there exists a connected two-fold covering of this group. The  covering  space is  called metaplectic group, and it is usually denoted by $Mp(2l,\mathbb{R}).$
Let us fix an element of the 
isomorphism class of all connected $2 : 1$ coverings of $Sp(2l,\mathbb{R})$ 
and denote it by $\lambda$. In particular, the mapping $\lambda: Mp(2l,\mathbb{R}) \to Sp(2l,\mathbb{R})$ is a Lie group
 representation.
 A  metaplectic structure on a symplectic manifold $(M^{2l},\omega)$ is a notion parallel to that one of a spin structure known from 
Riemannian geometry.
 In particular, one of its part is a principal $Mp(2l,\mathbb{R})$-bundle $\mathcal{Q}$ covering twice the bundle of 
symplectic repers $\mathcal{P}$ on $(M,\omega).$ This principal $Mp(2l,\mathbb{R})$-bundle is the mentioned metaplectic bundle and we will 
suppose this bundle to be chosen, and keep denote it by $\mathcal{Q}$ throughout this section.
  
As we have already said, the fields we shall be interested in, are certain exterior differential forms on $M^{2l}$ with values in the  {\it symplectic spinor bundle}. 
The symplectic spinor bundle is a  vector bundle over $M$ associated to 
the chosen principal $Mp(2l,\mathbb{R})$-bundle $\mathcal{Q}$ via an 'analytic derivate' of the Segal-Sahle-Weil representation. 
The Segal-Shale-Weil 
representation, denoted by $\tilde{L}$ in this text, is a faithful unitary representation of the metaplectic group
$Mp(2l,\mathbb{R})$ on the vector space $L^2(\mathbb{L})$ of  complex valued square Lebesgue integrable functions defined on a Lagrangian subspace $\mathbb{L}$ of the (standard) symplectic vector space 
$(\mathbb{R}^{2l},\omega_0).$ 
For technical reasons, we shall use the so called Casselman-Wallach globalization $V_{\infty}(HC(L^2(\mathbb{L})))$ of the underlying Harish-Chandra 
$(\mathfrak{g},\tilde{K})$-module of the Segal-Shale-Weil representation. (Here, $\mathfrak{g}$ is the Lie algebra of the symplectic group $G$ and $\tilde{K}$ is a maximal compact subgroup of the metaplectic group $\tilde{G}.$)
 We shall denote the resulting representation by  by $L$ and
call it the {\it metaplectic representation}. 
Thus $L: Mp(\mathbb{V},\omega_0) \to \mbox{Aut}({\bf S}),$ where ${\bf S}:=V_{\infty}(HC(L^2(\mathbb{L}))).$
We always consider the vector space $\bf S$ to be equipped with the action of the group $Mp(2l,\mathbb{R})$ and
call its {\it symplectic spinors}.  Let us mention that ${\bf S}$ decomposes into two irreducible $Mp(2l,\mathbb{R})$-submodules ${\bf S}_+$ and ${\bf S}_-,$ i.e., 
${\bf S} = {\bf S}_+ \oplus {\bf S}_-.$
In the following paragraph, we shall briefly explain why the  elements of $\bf S$ are called (symplectic) spinors. To our knowledge, the term symplectic spinor was first used by B. Kostant \cite{Kostant}.

It is known that the infinitesimal $\mathfrak{mp}(2l,\mathbb{R})$-module structure of the underlying Harish--Cha-ndra 
$(\mathfrak{g},\widetilde{K})$-module of ${\bf S}$ is equivalent to the space of polynomials $\mathbb{C}[x^1,\ldots, x^l]$ on which the Lie algebra $\mathfrak{mp}(2l,\mathbb{R})$ acts by the so called Dixmier representation
of the symplectic Lie algebra $\mathfrak{sp}(2l,\mathbb{R})\simeq \mathfrak{mp}(2l,\mathbb{R}).$ 
(The Dixmier  representation is an injective Lie algebra homomorphism of the Lie algebra $\mathfrak{mp}(2l,\mathbb{C})$ into the Lie algebra $\mbox{End}(\mathbb{C}[x^1,\ldots, x^l]).$)  Because $\mathbb{C}[x^1,\ldots,
x^l] \simeq \bigoplus_{k=0}^{\infty}S^{k}(\mathbb{L}^{\mathbb{C}})$ (the {\it symmetric  power} of 
the complexification of the Lagrangian subspace $\mathbb{L}$), we see that the situation is parallel ("super-symmetric")
to the case of $\mathfrak{so}(2l, \mathbb{C}),$ for which  spinors can be considered as elements of the {\it exterior power}
 of a maximal isotropic subspace in $\mathbb{C}^{2l}$   wr. to the bilinear form  
defining the Lie algebra  $\mathfrak{so}(2l,\mathbb{C}).$  Besides the super-symmetry, the fact that the metaplectic 
representation does not 
descend to a representation of the symplectic group gives a further justification for the use of the term 'spinor' in this case.
See Howe \cite{Howe} and Shale \cite{Shale} for related information. In the latter reference, one can find a use of  symplectic spinors, i.e., elements in $\bf S,$ in a quantization of boson fields.

  After we have introduced the Segal-Shale-Weil and metaplectic representations, let us describe the mentioned fields on which the operators 
from the symplectic twistor complexes  act more precisely. 
 The underlying algebraic structure of the symplectic spinor valued exterior differential forms 
is the vector space of  symplectic spinor valued exterior  forms, i.e., the vector space 
${\bf E}=\bigwedge^{\bullet}(\mathbb{R}^{2l})^*\otimes {\bf S}.$ We are considering the Grothendieck tensor product topology on $\bf E.$  Obviously, this vector space is equipped with the following 
tensor product representation $\rho$ of the metaplectic group $Mp(2l,\mathbb{R}).$  Thus, for $r=0,\ldots, 2l,$ $g\in Mp(2l,\mathbb{R})$ and $\alpha \otimes s
\in \bigwedge^r(\mathbb{R}^{2l})^*\otimes {\bf S},$ we set $\rho(g)(\alpha \otimes s):=\lambda(g)^{* \wedge r}\alpha \otimes L(g)s$
 and extend this prescription linearly. With this notation in mind, 
 the symplectic spinor valued exterior differential forms are sections of the vector bundle $\mathcal{E}$ associated to the chosen principal $Mp(2l,\mathbb{R})$-bundle $\mathcal{Q}$ via $\rho,$ i.e., $\mathcal{E}:=\mathcal{Q} \times_{\rho} {\bf E}.$

Now, we shall restrict our attention to the mentioned specific symplectic spinor valued exterior differential forms. 
For each $r=0,\ldots, 2l,$ there exists a distinguished irreducible submodule of $\bigwedge^{r}(\mathbb{R}^{2l})^* \otimes {\bf S}_{\pm}$ 
which we denote by ${\bf E}^{r}_{\pm}.$ Actually, the submodules ${\bf E}^r_{\pm}$ are the 
Cartan components of $\bigwedge^{r}(\mathbb{R}^{2l})^* \otimes {\bf S}_{\pm},$ i.e., their highest weight is the largest one of  the highest weights of all irreducible constituents of $\bigwedge^r(\mathbb{R}^{2l})^* \otimes {\bf S}_{\pm}$ wr. to the classical choices.   For $r=0,\ldots, 2l,$ we set ${\bf E}^r:={\bf E}^r_+ \oplus {\bf E}^r_-$ and $\mathcal{E}^r:=\mathcal{Q}\times_{\rho} {\bf E}^r.$
Further, let us denote the corresponding $Mp(2l,\mathbb{R})$-equivariant projection from $\bigwedge^{r} (\mathbb{R}^{2l})^* \otimes {\bf S}$ onto
${\bf E}^r$ by $p^r.$ We denote the lift of the projection 
 $p^r$ to the associated structures by the same symbol, i.e.,
 $p^r: \Gamma(M,\mathcal{Q}\times_{\rho} (\bigwedge^r \mathbb{R}^{2l} \otimes {\bf S})) \to \Gamma(M,\mathcal{E}^r).$ 


Now, we are in a position to define the subject of our investigation, namely the symplectic twistor complexes.
 Let us consider a Fedosov manifold  $(M, \omega, \nabla)$ and suppose 
that $(M,\omega)$ admits a metaplectic structure. Let $d^{\nabla^S}$ be the exterior covariant derivative  associated to $\nabla.$ 
For each $r=0,\ldots, 2l,$ let us
restrict the associated exterior covariant derivative $d^{\nabla^S}$ to $\Gamma(M, \mathcal{E}^{r})$ and compose the restriction with the projection 
$p^{r+1}.$ The resulting operator will be called  symplectic twistor operator and we
 will denote it by  $T_r.$ In this way, we obtain two sequences, namely 
$0  \longrightarrow
  \Gamma(M,\mathcal{E}^{0}) \overset{T_0}{\longrightarrow} 
  \Gamma(M,\mathcal{E}^{1}) \overset{T_{1}}{\longrightarrow}  
  \cdots  \overset{T_{l-1}}{\longrightarrow}
  \Gamma(M,\mathcal{E}^{l}) \longrightarrow 0 \mbox{   and}  
$ 
$0  \longrightarrow
  \Gamma(M,\mathcal{E}^{l}) \overset{T_l}{\longrightarrow} 
  \Gamma(M,\mathcal{E}^{l+1}) \overset{T_{l+1}}{\longrightarrow}  
  \cdots  \overset{T_{2l-1}}{\longrightarrow}
  \Gamma(M,\mathcal{E}^{2l}) \longrightarrow 0.  
$
It is known, see  Kr\'ysl \cite{SC}, that these   sequences are complexes provided the Fedosov manifold $(M^{2l},\omega,\nabla)$ 
is of Ricci type. These two complexes are the mentioned {\it symplectic twistor complexes}.
Let us notice, that we did not choose the full sequence of all symplectic spinor valued exterior differential forms and the exterior covariant derivative acting between them  because for a Ricci type Fedosov manifold, this sequence would not form a complex in general.

As we have mentioned, we shall prove that some parts of these two complexes are elliptic. 
To obtain these parts, one  should remove the last (the zero) term and the second last term from the first complex and the first (the zero) term from the second complex. The complexes obtained in this way will be called {\it left} and {\it right truncated symplectic twistor complex} according whether we have removed the  bundles from the first  symplectic twistor complex or from the second one, respectively.  
 Let us mention that by an elliptic complex, we mean a complex of differential operators such that its associated symbol sequence is an exact sequence of the sheaves in question. (We shall make this definition more precise in the text. See Schulze et al. \cite{Schulze} for details.)

Let us make some remarks on the methods we used to prove the ellipticity of the symplectic twistor complexes. We decided to use the so called Schur-Weyl-Howe correspondence, which we refer to as Howe correspondence for simplicity and which 
assigns to each representation of a Lie group $G$ the so called dual partner and certain representation of this partner.
In general, the Howe duality helps us to treat the representations of the group we started with. In the case of $GL(\mathbb{V})$ acting on $\mathbb{V}^{\otimes k}$ via the tensor product of its defining representation, the dual pair is the symmetric group $\mathfrak{S}_k$ on $k$ letters acting on $\mathbb{V}^{\otimes k}$ by permuting the positions of vectors constituting the appropriate $k$-tensor.
The presence of the  symmetric group, combinatorial in its nature, leads to several  combinatorial tools which are  well known in the representation theory of the general linear group $GL(\mathbb{V}),$ e.g., to  the concept of Young diagrams. The Howe type correspondence in our case, i.e., for the metaplectic group $Mp(2l,\mathbb{R})$ acting on the space $\bf E$ of symplectic spinor valued exterior forms, leads to the "smallest" simple super Lie algebra, namely to the ortho-symplectic super Lie algebra $\mathfrak{osp}(1|2)$ and certain representation on $\bf E.$  
We decided to use the Howe type correspondence mainly because the spaces ${\bf E}^r$ can be 
characterized via the mentioned representation  of $\mathfrak{osp}(1|2)$ easily. 
See R. Howe  \cite{Howe} for more information on the Howe type correspondence in general \cite{Howe}. 
Let us also mention that we have used the Cartan lemma on exterior differential forms in the proof of the ellipticity.

For other examples of elliptic complexes, we refer the reader, e.g., to
Hotta \cite{Hotta}, Schmid \cite{Schmid} or \cite{Schmid2} and Stein and Weiss \cite{SW}. Let us mention
 that the proofs of the ellipticity of the deRham and Dolbeault complexes are also based on a use of the 
Cartan lemma only. Roughly said, this is mainly because of the relatively simple representation theory 
of the orthogonal (deRham case) and unitary groups (Dolbeault case) on exterior forms and
their complexification, respectively. See, e.g., Wells \cite{Wells}.

  In the second section, we recall some basic facts from symplectic linear algebra, mention facts on symplectic spinors 
and symplectic spinor valued exterior forms and its decomposition into irreducible submodules.
In the third chapter, basic facts on Fedosov manifolds and their curvature are mentioned and the symplectic twistor complexes are introduced. In the fourth chapter, the symbol sequence of the symplectic twistor complexes is computed,
several lemmas of technical character are derived and finally, the ellipticity of the truncated symplectic twistor complexes is proved.

\section{Symplectic spinor valued forms}

In this section, we set a notation, recall some facts from symplectic
 linear algebra, give a definition the metaplectic group  and introduce the basic object of our study, 
namely the space of symplectic spinor valued exterior forms. In the whole text the Einstein summation convention is used, 
 not mentioning it explicitly. 

\subsection{ {\it Symplectic group and its action on exterior forms }}

In order to set the notation, let us start recalling some simple results from  symplectic linear algebra.
Let $(\mathbb{V},\omega_0)$ be a real symplectic vector space of dimension $2l,$ $l \geq 1.$  Let us choose two Lagrangian 
subspaces $\mathbb{L}$ and $\mathbb{L}',$  such that $\mathbb{V} \simeq \mathbb{L} \oplus \mathbb{L}'$ \footnote{Let us recall 
that by Lagrangian, we mean maximal isotropic wr. to $\omega_0$}. It is easy to check that
$\dim \mathbb{L} = \dim \mathbb{L}'=l.$
Let us choose an adapted symplectic basis $\{e_i\}_{i=1}^{2l}$ of $(\mathbb{V} \simeq \mathbb{L} \oplus \mathbb{L}', \omega_0),$
i.e., $\{e_i\}_{i=1}^{2l}$ is a symplectic basis of $(\mathbb{V},\omega_0)$ and $\{e_{i}\}_{i=1}^l \subseteq \mathbb{L}$ and 
$\{e_i\}_{i=l+1}^{2l} \subseteq \mathbb{L}'.$ 
The basis dual to the basis $\{e_i\}_{i=1}^{2l}$ will be denoted by $\{\epsilon^i\}_{i=1}^{2l},$ i.e., for $i,j=1,\ldots, 2l$ we 
have $\epsilon^j(e_i)=
\iota_{e_i}\epsilon^j=\delta^{j}_i,$ where $\iota_v\alpha$ for an element $v\in \mathbb{V}$ and an exterior form $\alpha \in 
\bigwedge^{\bullet}\mathbb{V}^*,$
denotes the contraction of the form $\alpha$ by the vector $v.$
 Further for $i, j = 1, \ldots, 2l,$ we set
$\omega_{ij}:= \omega_0(e_i,e_j)$ and define $\omega^{ij},$ $i,j=1,\ldots, 2l,$ by the equation $\omega_{ij}\omega^{kj} = \delta_i^k$ for all
$i,k =1, \ldots, 2l.$ (Here the summation convention was used.) Let us remark that not
 only $\omega_{ij}= -\omega_{ji},$ but also $\omega^{ij}=-\omega^{ji}$ for $i,j =1, \ldots, 2l.$

As in the Riemannian case, we would like to rise and lower indices of tensor coordinates.
 In the symplectic case, one should be more careful because of the anti-symmetry of $\omega_0.$ 
For coordinates ${K_{ab\ldots c\ldots d}}^{rs \ldots t \ldots u}$ of a tensor $K$ over $\mathbb{V},$  we denote
the expression $\omega^{ic}{K_{ab\ldots c \ldots d}}^{rs \ldots t}$ by 
${{{K_{ab \ldots}}^{i}}_{\ldots d}}^{rs \ldots t}$ and 
${K_{ab\ldots c}}^{rs \ldots t \ldots u}\omega_{ti}$ by ${{{K_{ab \ldots c}}^{rs\ldots}}_{i}}^{\ldots u}$ 
and similarly for other types of tensors and also in the geometric setting when we will be considering tensor
 fields over a symplectic manifold $(M,\omega)$.   
Further, one can also define an isomorphism $\sharp: \mathbb{V}^* \to \mathbb{V} \mbox{,  } \alpha \mapsto \alpha^{\sharp},$  by
the formula $$\alpha(w)= \omega_0(\alpha^{\sharp},w) \mbox{ for each }\alpha \in \mathbb{V}^* \mbox{ and }w \in \mathbb{V}.$$ 
For $\alpha = \alpha_i\epsilon^i$ and
$j=1,\ldots, 2l,$ we get $\alpha_j= \alpha(e_j)=\omega_0((\alpha^{\sharp})^ie_i,e_j)=\omega_{ij}(\alpha^{\sharp})^i=(\alpha^{\sharp})_{j}$ 
which implies  $\alpha^{\sharp} = (\alpha^{\sharp})^{i} e_i= \alpha^i e_i.$
Thus, we see that  the isomorphism $\sharp$ is realized by   rising of indices via the form $\omega_0.$

Now, let us introduce the groups we will be dealing with. Let us denote the symplectic group of 
$(\mathbb{V},\omega_0)$ by $G,$ i.e., $G:=Sp(\mathbb{V},\omega_0) \simeq Sp(2l,\mathbb{R}).$
Because the homotopy group of $G=Sp(\mathbb{V},\omega_0)$ is $\mathbb{Z},$ there exists a connected $2:1$ 
(necessarily non-universal) covering of 
$G$ by a Lie group $\tilde{G},$ the so called metaplectic group; here denoted by $\tilde{G} := Mp(\mathbb{V},\omega_0) 
\simeq Mp(2l,\mathbb{R}).$ 
Denote by  $\lambda: \tilde{G} \to G$  the mentioned two-fold covering. 

\subsection{ {\it Segal-Shale-Weil representation and symplectic spinor valued forms}  }

  The Segal-Shale-Weil representation is a distinguished representation of the metaplectic group $\tilde{G}=Mp(\mathbb{V},\omega_0).$\footnote{The names oscillator and metaplectic are also used in the literature. See, e.g., Howe \cite{Howe}.} As we have mentioned in the Introduction, this representation is unitary, faithful and does not descend to a representation of the symplectic group.  
Its underlying vector space is the vector space of 
complex valued square Lebesgue integrable functions $L^2(\mathbb{L})$ defined on the Lagrangian subspace $\mathbb{L}.$ 
To set a notation, let us denote the Segal-Shale-Weil representation by $\tilde{L},$ i.e., 
$$\tilde{L}: \tilde{G}\to \mathcal{U}(L^2(\mathbb{L})),$$ where $\mathcal{U}(H)$ denotes the unitary group of a Hilbert space $H.$ Let 
us set ${\bf S}:=V_{\infty}(HC(L^2(\mathbb{L}))),$ where  $V_{\infty}$ is the Casselma-Wallach globalization functor and $HC$ be the forgetful Harish-Chandra functor.
We shall denote the resulting Casselman-Wallach globalization of the Segal-Shale-Weil representation by $L$ and call it the 
{\it metaplectic representation} and the elements of $\bf S$ {\it symplectic spinors}.  
It is well known that ${\bf S}$ splits into two irreducible $Mp(\mathbb{V},\omega_0)$-submodules ${\bf S}_+$ and ${\bf S}_-.$
Thus, we have ${\bf S} = {\bf S}_+ \oplus {\bf S}_-.$ See Weil \cite{Weil} and Kashiwara, Vergne \cite{KV} for more detailed information 
on the Segal-Shale-Weil representation and Casselman \cite{CW} on this type of globalization.

  Now, we may define the so called {\it symplectic Clifford multiplication} $\cdot  : \mathbb{V} \times {\bf S} \to {\bf S}.$
For $s\in {\bf S},$ $x=x^je_j \in \mathbb{L},$ $x^j \in \mathbb{R},$ $j=1,\ldots, 2l$ and $i=1,\ldots, l,$ let us set
\begin{eqnarray*}
e_i.s (x) &:=& \imath x^i s (x) \mbox{ and }\\
e_{i+l}. s(x) &:=& \frac{\partial s}{\partial x^{i}}(x).
\end{eqnarray*}
In physics, this mapping (up to a constant multiple) is usually called the canonical quantization.

  Let us remark that the definition is correct  because of an interpretation of the Casselman-Wallach (also called smooth) globalization 
See, e.g., Vogan \cite{Vogan} for details on this interpretation.

 For each $v,w \in \mathbb{V}$ and $s\in {\bf S},$ one can easily derive the following commutation relation
\begin{eqnarray}
v.w.s - w.v.s = -\imath \omega_0(v,w)s. \label{cr}
\end{eqnarray}
(For a proof, see, e.g., Habermann, Habermann \cite{HH}.) We shall use this relation repeatedly and almost always without mentioning its use.  

Now, we  prove that the symplectic Clifford multiplication by a fixed non-zero vector $v\in \mathbb{V}$
is injective as a mapping from ${\bf S} \to {\bf S}.$
We shall use the equivariance of the symplectic Clifford multiplication, i.e., the fact $L(g)(v.s) = [\lambda(g)v].L(g)s$ which hold for each $g  \in \tilde{G},$ $v\in \mathbb{V}$ and $s \in {\bf S}.$ (See Habermann, Habermann \cite{HH}.) Thus, let us suppose that $s \in {\bf S}$ and $0 \neq v \in \mathbb{V}$ are given such that
$v.s = 0.$ Because the action of the symplectic group $G$ on $\mathbb{V}-\{0\}$ is transitive and $\lambda$ is a covering, there exists an element $g \in \tilde{G}$ such that $\lambda(g)v = e_1.$ Applying $L(g)$ on the equation 
$v.s = 0,$ we get $L(g)(v.s) = 0.$ Using the above mentioned equivariance of the symplectic Clifford multiplication, we get $0= L(g)(v.s) = [\lambda(g)v].(L(g)s) = e_1.(L(g)s).$ Denoting $L(g)s =: \psi$ and using the definition of the symplectic Clifford multiplication, we get $x^1\psi = 0,$ which implies $\psi(x) = 0$ for each $x=(x^1,\ldots, x^l)  \in \mathbb{L}$ such that $x^1\neq 0.$
By continuity of $\psi,$ we get
$\psi = 0.$ Because $L$ is a group representation,  we get $s=0$ from $0=\psi = L(g)s,$ i.e., the injectivity of the symplectic Clifford multiplication. 

  Having defined the Segal-Shale-Weil representation and the symplectic Clifford multiplication, we shall introduce (the algebraic and analytic version 
of) the basic geometric structure we are be interested in. Namely, we  introduce  the space  ${\bf E}$ of symplectic spinor valued exterior forms, i.e., the space
$\bigwedge^{\bullet}\mathbb{V}^*\otimes \bf{S}.$  We shall consider this space to be equipped with the Grothendieck tensor product topology, cf.,
 e.g., Tr\'eves \cite{Treves}.

The metaplectic group $\tilde{G}:=Mp(\mathbb{V},\omega_0)$ acts on ${\bf E}$ by the representation
$$\rho: \tilde{G} \to \mbox{Aut}({\bf E}) \mbox{ defined by the formula}$$ 
$$\rho(g)(\alpha\otimes s):=(\lambda(g)^*)^{\wedge r}\alpha \otimes L(g)s,$$ 
where $\alpha\in \bigwedge^r\mathbb{V}^*,$ $s\in {\bf S},$ $r=0,\ldots, 2l,$ and it is
extended by linearity also to non-homogeneous elements.

For a vector $v\in \mathbb{V}$ and a homogeneous symplectic spinor valued exterior form $\psi:=\alpha \otimes s,$ we set
$\iota_v \psi := \iota_v \alpha \otimes s$ and $v.\psi:= \alpha \otimes v.s$ and extend the definition by linearity to 
non-homogeneous elements. 

Now, we shall  describe the decomposition of $\bf E$ into irreducible $Mp(\mathbb{V},\omega_0)$-submodules.
 For $i=0,\ldots, l,$ let us set $m_i:=i,$ and for $i=l+1,\ldots 2l,$ we set $m_i:=2l-i,$ 
and define the   set $\Xi$ of pairs of non-negative integers
$$\Xi:=\{(i,j) \in \mathbb{N}_0\times \mathbb{N}_0| i=0,\ldots, 2l, j=0,\ldots, m_i\}.$$
One can say the set $\Xi$  has a shape of a triangle if visualized in a $2$-plane. (See the Figure 1. bellow). 
We use this set for parameterizing the irreducible submodules of  $\bf E.$

In Kr\'ysl \cite{KryslJOLT2} for each $(i,j) \in \Xi,$ two irreducible $\tilde{G}$-modules ${\bf E}^{ij}_{\pm}$ were
uniquely defined via the highest weights of their underlying Ha\-rish-Chan\-dra module and the fact that they are irreducible submodules of 
$\bigwedge^i\mathbb{V}^*\otimes {\bf S}_{\pm}.$  
For our convenience for each
 $(i,j) \in \mathbb{Z}\times \mathbb{Z} \setminus \Xi,$ let us put ${\bf E}^{ij}_{\pm}:=0,$ and for each 
$(i,j) \in \mathbb{Z}\times \mathbb{Z},$  we define ${\bf E}^{ij}:={\bf E}_{+}^{ij} \oplus {\bf E}_{-}^{ij}.$
In general, all objects equipped with one index or a tuple of indices are supposed to be zero if the index is out of the range $1,...,2l$ or the tuple is out of the set $\Xi,$ respectively.

In the following theorem, the decomposition of ${\bf E}$ into irreducible $\tilde{G}$-sub--modules is described.

{\bf Theorem 1:} For $r=0,\ldots, 2l,$ the following decomposition into irreducible 
$\tilde{G}$-modules
$$\bigwedge^r \mathbb{V}^{*} \otimes {\bf S}_{\pm} \simeq \bigoplus_{j, (r,j)\in\Xi}  {\bf E}^{rj}_{\pm} \quad \mbox{ holds}.$$

{\it Proof.} See Kr\'ysl \cite{KryslJOLT2}. $\Box$

Whereas it is not necessary to provide the reader with a prescription for the highest weights of the $\tilde{G}$-modules ${\bf E}^{ij}_{\pm},$ 
the following remark on the multiplicity structure of the module ${\bf E}$ will be crucial.  This remark follows from the prescriptions for the highest weights of the infinitesimal structure of the underlying Harish--Chandra modules of
${\bf E}_{\pm}^{ij}$ in Kr\'{y}sl \cite{KryslSVF} easily.

{\bf Remark:} \begin{itemize}
\item[1.]   For any $(r,j),(r,k) \in \Xi$ such that $j \neq k,$ we have
$${\bf E}^{rj}_{\pm} \not \simeq {\bf E}^{rk}_{\pm}$$ (any combination of $\pm$ 
at both sides of the preceding
relation is allowed).
Thus in particular, $\bigwedge^r \mathbb{V}^* \otimes {\bf S}$ is multiplicity-free for each  $r = 0, \ldots, 2l.$
\item[2.] Moreover, it is known that ${\bf E}^{rj}_{\pm} \simeq {\bf E}^{sj}_{\mp}$ for each $(r,j), (s,j) \in \Xi.$ One cannot change the order of $+$ and $-$   at precisely one side of the preceding isomorphism without changing its trueness.
\item[3.] From the preceding two items, one gets immediately that there are no submodules of $\bigwedge^i \mathbb{V}^* \otimes {\bf S}$ isomorphic to 
${\bf E}^{i+1,i+1}_{\pm}$ for each $i=0,\ldots, l-1.$
\end{itemize}

  In the next figure (Figure 1.), one can see the decomposition structure of $\bigwedge^{\bullet}\mathbb{V}^{*}\otimes {\bf S}_{\pm}$ in the case of $l=3.$ For $i=0,\ldots, 6,$ the $i^{th}$ column constitutes of the irreducible modules in which the ${\bf S}_{\pm}$-valued exterior 
forms of form-degree $i$ decompose.
 
$$\xymatrix{
{\bf E}_{\pm}^{00}  &{\bf E}^{10}_{\pm}  &{\bf E}^{20}_{\pm}  &{\bf E}^{30}_{\pm}  &{\bf E}^{40}_{\pm}  &{\bf E}^{50}_{\pm}  &{\bf E}^{60}_{\pm} \\
&{\bf E}^{11}_{\pm}   &{\bf E}^{21}_{\pm} &  {\bf E}^{31}_{\pm}   & {\bf E}^{41}_{\pm} &{\bf E}^{51}_{\pm}  & \\
& & {\bf E}^{22}_{\pm}   & {\bf E}^{32}_{\pm}   & {\bf E}^{42}_{\pm} &&\\
&&& {\bf E}^{33}_{\pm} &&&}
$$

\centerline{Figure 1. Decomposition of $\bigwedge^{\bullet}\mathbb{V}^* \otimes {\bf S}_{\pm}$ for $l=3.$}

In the next theorem, the decomposition of $\mathbb{V}^* \otimes {\bf E}^{ij},$ $(i,j) \in \Xi,$ is described. 
Let us remind the reader that due to our convention ${\bf E}^{ij} = 0$ for $(i,j) \in \mathbb{Z}\times \mathbb{Z} \setminus \Xi.$  

{\bf Theorem 2:} For $(i,j) \in \Xi,$ we have 
$$(\mathbb{V}^*\otimes {\bf E}^{ij}) \cap (\bigwedge^{i+1}\mathbb{V}^*\otimes {\bf S})
 \simeq {\bf E}^{i+1,j-1}\oplus {\bf E}^{i+1,j} \oplus {\bf E}^{i+1,j+1}.$$
{\it Proof.} See  the proof of the Theorem 4 in Kr\'ysl \cite{KryslSVF}. $\Box$

{\bf Remark:} In the mentioned Theorem 4,  a superset of the image of an exterior covariant
 derivative was determined. The problem of determining  this superset was reduced 
to the question of determining $(\mathbb{V}^*\otimes {\bf E}^{ij}) 
\cap (\bigwedge^{i+1}\mathbb{V}^*\otimes {\bf S})$ immediately. 
The intersection computed there is identical to the right-hand side of the isomorphism in the Theorem 2 of this paper.

\subsection{ {\it Operators related to a Howe type correspondence.}}

  In this section, we will introduce five 
 operators acting on the space ${\bf E}$ of symplectic spinor valued exterior forms.
 These operators are related to the so called Howe type correspondence for the  metaplectic group $Mp(\mathbb{V},\omega_0)$ 
acting on ${\bf E}$ via the representation $\rho$ introduced above.
  
For   $ r=0,\ldots, 2l$ and $\alpha \otimes s \in \bigwedge^r
\mathbb{V}^*\otimes {\bf S},$  we set
$$F^+:\bigwedge ^{r}\mathbb{V}^*\otimes {\bf S} \to
\bigwedge^{r+1}\mathbb{V}^*\otimes {\bf S},\, F^+(\alpha \otimes
s):=\frac{\imath}{2}\sum_{i=1}^{2l}\epsilon^i\wedge \alpha \otimes e_i.s,$$
$$F^-:\bigwedge ^{r}\mathbb{V}^* \otimes {\bf S} \to \bigwedge ^{r-1}\mathbb{V}^*\otimes {\bf S}, \,
F^-(\alpha \otimes s):=\frac{1}{2}\sum_{i=1}^{2l}\omega^{ij}\iota_{e_i}\alpha \otimes
e_j.s$$
and extend them linearly.

Next, we shall define the operators $H, E^+$ and $E^-.$
For $r=0,\ldots, 2l,$ we set
$$H:\bigwedge ^{r}\mathbb{V}^*\otimes {\bf S} \to \bigwedge ^{r}\mathbb{V}^*\otimes {\bf S},\, H:=2\{F^+,F^-\} \, \mbox{ and}$$
$$E^{\pm}:\bigwedge ^{r}\mathbb{V}^*\otimes {\bf S} \to \bigwedge ^{r\pm 2}\mathbb{V}^*\otimes {\bf S}, \,
E^{\pm}:=\pm 2 \{F^{\pm},F^{\pm}\},$$ where $\{,\}$ denotes the anti-commutator in the associative algebra $\mbox{End}({\bf E}).$

In the next lemma, we sum-up some known facts on the operators $F^{\pm},$ $E^{\pm}$ and $H$ which we shall need in the proof of 
the ellipticity of the truncated symplectic twistor complexes.

{\bf Lemma 3}: \begin{itemize}
\item[1.] The operators $F^{\pm},$ $E^{\pm}$ and $H$ are $\tilde{G}$-equivariant.
\item[2.] For $i=0,\ldots, l,$ the operator $F^-_{|{\bf E}^{im_i}} = 0.$
\item[3.] For $\alpha \otimes s \in \bigwedge^{\bullet}\mathbb{V}^* \otimes {\bf S},$ we have
\begin{eqnarray}
E^-(\alpha \otimes s)&=&\frac{\imath}{2}\omega^{ij}\iota_{e_i}\iota_{e_j}\alpha \otimes s. \label{ecko}
\end{eqnarray} 
\item[4.] The associative algebra\\
 $\mbox{End}_{\tilde{G}}({\bf E}):=\{A \in \mbox{End}({\bf E}) | A\rho(g) = \rho(g)A \mbox{ for all } g \in \tilde{G}\}$
is generated by $F^+$ and $F^-$ and the projections $p_{\pm},$ where $p_{\pm}:{\bf S} \to {\bf S}_{\pm}$ are the $\tilde{G}$-equiavariant projections.
\end{itemize}
{\it Proof.} See Kr\'ysl \cite{KryslJOLT2}. $\Box$

According to the third item of the preceding lemma, we see that the operator $E^-$ acts on the form-part of a symplectic spinor
 valued exterior form only. Because of that, we will make no difference between 
$E^-(\alpha \otimes s)$ and $E^-(\alpha) \otimes s$ and we will   write $E^- \alpha \otimes s$ instead of 
any of the previous two expressions briefly.

Now, let us prove a lemma of a technical character. We shall use it when we will be proving the ellipticity of the truncated symplectic twistor complexes.

{\bf Lemma 4:}  For each $r=0,\ldots, 2l,$ $(i,j)\in \Xi,$ $v\in \mathbb{V}$ and $\alpha \otimes s \in \bigwedge^r \mathbb{V}^* \otimes {\bf S}$ the following relations 
\begin{eqnarray}
[E^+,E^-]=H, &&  [E^-, F^+] = -F^- \label{ef},\\
H(\alpha\otimes s)&=& \frac{1}{2}(r-l)\alpha \otimes s \label{h},\\
F^-F^+_{|{\bf E}^{ij}} &=& \frac{1}{4}(\frac{1+i-j}{2})\mbox{Id}_{{\bf E}^{ij}} \mbox{ for } i+j \mbox{ odd} \label{ffo},\\
F^-F^+_{|{\bf E}^{ij}} &=& \frac{1}{4}(\frac{i+j}{2}-l)\mbox{Id}_{{\bf E}^{ij}} \mbox{ for } i+j \mbox{ even} \label{ffec},\\ 
\{F^+,\iota_v \}(\alpha \otimes s)=\frac{\imath}{2} \alpha \otimes v.s &\mbox{ and}& [F^-,  v. ](\alpha \otimes s)=\frac{\imath}{2}\iota_v\alpha \otimes s  \label{fc}
\end{eqnarray}
hold.

{\it Proof.} See Kr\'ysl \cite{KryslJOLT2} for a proof of the relations in the rows (\ref{ef}), (\ref{h}), (\ref{ffo}) and (\ref{ffec}).
Let us suppose we are given an element $v=v^ie_i\in \mathbb{V},$ $v^i \in \mathbb{R},$ $i=1,\ldots, 2l,$ and a homogeneous 
element $\alpha \otimes s \in \bigwedge^j\mathbb{V}^*\otimes {\bf S},$ $j=0,\ldots, 2l$.
First, let us prove the first relation in the row (\ref{fc}). Using the definition of $F^+$, we may write
$\{F^+, \iota_v\}(\alpha \otimes s)= F^+ (\iota_v\alpha \otimes s) + \frac{\imath}{2}\iota_v (\epsilon^i \wedge \alpha \otimes e_i.s)=
\frac{\imath}{2}[\epsilon^i \wedge \iota_v \alpha \otimes e_i.s + v^i \alpha \otimes e_i.s - \epsilon^i \wedge \iota_v \alpha \otimes e_i.s]
=\frac{\imath}{2}\alpha\otimes v.s.$ Thus, the first relation of (\ref{fc}) follows by linearity. 
Now, let us prove the second relation at the row (\ref{fc}). Using the definition of $F^-$ and 
the commutation relation (\ref{cr}), we get
$F^- (\alpha \otimes v.s) = \frac{1}{2}(\omega^{ij}\iota_{e_i}\alpha \otimes e_j.v.s)= 
\frac{1}{2}\omega^{ij}\iota_{e_i}\alpha \otimes (v.e_j.s - \imath \omega_0(e_j,v)s)= v. F^-(\alpha \otimes s) 
+ \frac{\imath}{2}\omega^{ij}\iota_{e_i}\alpha \otimes v_j s= v. F^{-}(\alpha \otimes s) + \frac{\imath}{2}\iota_v \alpha \otimes s.$ 
Thus, the second relation at the row (\ref{fc}) is proved. 
$\Box$

{\bf Remark:} The operators $F^{\pm}, E^{\pm}$ and $H$ satisfy the commutation and anti-commutation relations identical to that one 
which are satisfied 
by the usual generators of the ortho-symplectic super Lie algebra $\mathfrak{osp}(1|2).$  One can say that the super-Lie algebra 
$\mathfrak{osp}(1|2)$ is a Howe dual partner to the metaplectic group $Mp(\mathbb{V},\omega_0)$ acting on $\bf E$ via $\rho.$
See the Introduction and the references therein for more on the Howe type correspondence.




{




\section{Symplectic twistor complex}

In this section, we define the notion of a Fedosov manifold, recall some informations on its curvature, introduce the
symplectic analogue of a spin structure, namely the metaplectic structure) and define the symplectic twistor complex.

  Let $(M,\omega)$ be a symplectic manifold. Let us consider an affine torsion-free symplectic connection
$\nabla$ on $(M,\omega)$ and denote the induced connection on $\Gamma(M,\bigwedge^2 T^*M)$ also by $\nabla$.
Let us recall that by  torsion-free and symplectic, we mean 
$T(X,Y):=\nabla_XY-\nabla_YX -[X,Y]=0$ for all $X,Y \in \mathfrak{X}(M)$ and $\nabla \omega = 0$ 
 These connections are usually called Fedosov connections. See the Introduction and the references therein for an application of these connections. 
The curvature tensor $R^{\nabla}$ of a Fedosov connection   is defined in the classical way, i.e., formally by the same formula as in the Riemannian geometry.
It is known, see, e.g., Vaisman \cite{Vaisman}, that 
$R^{\nabla}$ splits into two  parts, namely into the symplectic Ricci and Weyl curvature tensors fields $\widetilde{\sigma}^{\nabla}$ 
and $W^{\nabla}.$ Let us display the definitions of these two curvature parts. In what follows whenever indices are not quantified 
and we are not summing over them, we suppose the appropriate predicate holds for all these indices from the set  $\{1,\ldots, 2l\}.$ For each $m\in M,$ we choose an open neighborhood  $U \ni m$ such that there exists an adapted symplectic frame $\{e_i\}_{i=1}^{2l}$ on it. Coordinates of a tensor field are always meant with respect to the chosen frame. 
For the coordinates of the symplectic Ricci and symplectic Weyl tensor fileds, we have (see Vaisman \cite{Vaisman})
\begin{eqnarray*}
\sigma_{ij} &:=&{R^{k}}_{ikj},\\
2(l+1)\widetilde{\sigma}^{\nabla}_{ijkl}&:=&\omega_{il}\sigma_{jk}-\omega_{ik}\sigma_{jl}+\omega_{jl}\sigma_{ik}
-\omega_{jk}\sigma_{il}+2\sigma_{ij}\omega_{kl} \mbox{ and}\\
 W^{\nabla}&:=&R^{\nabla} - \widetilde{\sigma}^{\nabla}.
\end{eqnarray*}
Let us call a Fedosov manifold of {\it Ricci type}, if $W^{\nabla}=0.$
 
{\bf Remark:}
Because the Ricci curvature tensor field is symmetric (see Vaisman \cite{Vaisman}), a 
possible symplectic scalar curvature does not exist because $\sigma^{ij}\omega_{ij}$ is zero. 

{\bf Example:}
It is easy to see that each Riemann surface equipped with its volume form as the symplectic form and with the Levi-Civita connection is a Fedosov manifold of Ricci type. Further 
for any $l\geq 1,$ the Fedosov manifold $(\mathbb{CP}^l,\omega_{FS},\nabla)$ is also a  Fedosov manifold of Ricci type. Here, 
$\omega_{FS}$ is the K\"{a}hler form associated to the Fubini-Study metric and to the standard complex 
structure on $\mathbb{CP}^l$, and $\nabla$ is the Levi-Civita connection associated to the Fubini-Study metric.

Now, let us  introduce the metaplectic structure the definition of which we have sketched briefly in the Introduction.
The metaplectic structure is a 
symplectic analogue of the notion of a spin structure in  the Riemannian geometry.
For a symplectic manifold $(M^{2l}, \omega)$  of dimension $2l,$
let us denote the bundle of symplectic reperes in $TM$ by
$\mathcal{P}$ and  the foot-point projection of $\mathcal{P}$ onto
$M$ by $p.$ Thus $(p:\mathcal{P}\to M, G),$ where $G\simeq
Sp(2l,\mathbb{R}),$ is a principal $G$-bundle over $M$. As in
the subsection 2.1., let $\lambda: \tilde{G}\to G$ be a member
of the isomorphism class of the non-trivial two-fold coverings of
the symplectic group $G.$ In particular, $\tilde{G}\simeq
Mp(2l,\mathbb{R}).$ Further, let us consider a principal
$\tilde{G}$-bundle $(q:\mathcal{Q}\to M, \tilde{G})$ over the chosen 
symplectic manifold $(M,\omega).$ We call the pair
$(\mathcal{Q},\Lambda)$   metaplectic structure if  $\Lambda:
\mathcal{Q} \to \mathcal{P}$ is a surjective bundle homomorphism
on $M$ and if the following diagram,
$$\begin{xy}\xymatrix{
\mathcal{Q} \times \tilde{G} \ar[dd]^{\Lambda\times \lambda} \ar[r]&   \mathcal{Q} \ar[dd]^{\Lambda} \ar[dr]^{q} &\\
                                                            & &M\\
\mathcal{P} \times G \ar[r]   & \mathcal{P} \ar[ur]_{p} }\end{xy}$$
with the
horizontal arrows being respective actions of the displayed groups, commutes.
See, e.g.,  Habermann,  Habermann \cite{HH} and Kostant \cite{Kostant} for
details on metaplectic structures. Let us only remark, that typical examples of symplectic
 manifolds admitting a metaplectic structure are cotangent bundles of orientable manifolds (phase spaces), Calabi-Yau manifolds and 
complex projective spaces $\mathbb{CP}^{2k+1}$, $k \in \mathbb{N}_0.$

Let us denote the  vector bundle  associated to the introduced principal $\tilde{G}$-bundle
$(q:\mathcal{Q}\to M,\tilde{G})$ via the metaplectic representation $L$ by $\mathcal{S}.$
Thus, we have $\mathcal{S}=\mathcal{Q}\times_{L}{\bf S}.$ We shall call
this associated vector bundle {\it symplectic spinor bundle}.  The 
sections $\phi \in \Gamma(M,\mathcal{S})$ will be called  {\it symplectic spinor fields}.
Let us put $\mathcal{E}:=\mathcal{Q}\times_{\rho}{\bf E}.$  For 
$r=0,  \ldots, 2l$ we define $\mathcal{E}^r:=\mathcal{Q}\times_{\rho} {\bf E}^r,$ where ${\bf E}^r$ abbreviates ${\bf E}^{rm_r}.$ 
The sections $\Gamma(M,\mathcal{E})$ will be called symplectic spinor valued exterior differential forms. (For the notion of 
differentiability of sections with values in Fr\'{e}chet bundles we use the metrizability of the fibers.) 


Because the operators $E^{\pm}, F^{\pm}$ and $H$  are 
$\tilde{G}$-equivariant (see the Lemma 3 item 1), they lift to operators
acting on sections of the corresponding associated vector bundles. The same is true about the projections 
$p^{ij},$ $(i,j) \in \mathbb{Z}\times \mathbb{Z},$ and $p_{\pm}$ because they are projections onto $\tilde{G}$-submodules.
We shall use the same symbols as for the defined operators as for
their "lifts" to the associated vector bundle structure.

 Now, we shall  make a use of the Fedosov connection.  
The Fedosov connection $\nabla$ determines the associated principal bundle connection $Z$
on the principal bundle $(p:\mathcal{P}\to M, G).$ 
This connection lifts to a principal bundle connection on  the principal bundle
$(q:\mathcal{Q}\to M, \tilde{G})$ and defines the associated covariant derivative 
on the symplectic bundle $\mathcal{S},$ which we shall denote by $\nabla^S,$ and call it the 
 symplectic spinor covariant derivative. See Habermann, Habermann \cite{HH} for details. 
The symplectic spinor covariant derivative $\nabla^S$
 induces the exterior symplectic spinor derivative $d^{\nabla^S}$ acting on $\Gamma(M,\mathcal{E}).$
For $r=0,\ldots, 2l,$ we have $d^{\nabla^S}: \Gamma(M,\mathcal{Q}\times_{\rho}(\bigwedge^r\mathbb{V}^* \otimes {\bf S}) \to \Gamma(M,\mathcal{Q}\times_{\rho} (\bigwedge^{r+1}\mathbb{V}^*\otimes {\bf S})).$  

Now, we can define the symplectic twistor operators.
For $r=0,\ldots, 2l,$ we set
$$T_r:\Gamma(M,\mathcal{E}^{r})\to \Gamma(M,\mathcal{E}^{r+1}), \quad \mbox{  } T_r:=
p^{r+1, m_{r+1}}d^{\nabla^S}_{|\Gamma(M,\mathcal{E}^{r})}$$ and call these operators  {\it symplectic twistor operators.} 
Informally, one can say that the operators are going on the edge of the triangle at the Figure 1. 
Let us notice that $F^-(\nabla^S-T_0)$ is, up to a nonzero scalar multiple, the so called symplectic Dirac operator 
introduced by K. Habermann. See, e.g., Habermann, Habermann \cite{HH}.

In the next theorem, we state the result on the two introduced symplectic twistor sequences. 

{\bf Theorem 5:} Let $l\geq 2$ and $(M^{2l},\omega,\nabla)$ be a Fedosov manifold of Ricci type admitting a metaplectic structure.
Then
$$0  \longrightarrow
  \Gamma(M,\mathcal{E}^{00}) \overset{T_0}{\longrightarrow} 
  \Gamma(M,\mathcal{E}^{11}) \overset{T_{1}}{\longrightarrow}  
  \cdots  \overset{T_{l-1}}{\longrightarrow}
  \Gamma(M,\mathcal{E}^{ll}) \longrightarrow 0 \mbox{   and}  
$$ 
$$0  \longrightarrow
  \Gamma(M,\mathcal{E}^{ll}) \overset{T_l}{\longrightarrow} 
  \Gamma(M,\mathcal{E}^{l+1,l+1}) \overset{T_{l+1}}{\longrightarrow}  
  \cdots  \overset{T_{2l-1}}{\longrightarrow}
  \Gamma(M,\mathcal{E}^{2l,2l}) \longrightarrow 0  
$$ 
are complexes.

{\it Proof.} See Kr\'ysl \cite{SC}. $\Box$






\section{Ellipticity of the symplectic twistor complex}

In this section, we prove the ellipticity of the truncated symplectic twistor complexes. Let us recall that by an elliptic complex of differential operators we mean a complex of differential operators such
 that the associated complex of symbols of the differential
 operators is an exact sequence 
of sheaves. (Let us notice that given a smooth vector bundle $\mathcal{V}$ over a smooth manifold $M,$ then we are always  considering the 
associated sheaf of smooth sections of this bundle.)
See, e.g., Wells \cite{Wells} and Schulze et al. \cite{Schulze} for more details on ellipticity.

We start with a simple lemma in which the symbol of the exterior covariant symplectic spinor derivative associated to a Fedosov manifold admitting a metaplectic structure is computed.

{\bf Lemma 6:} Let $(M,\omega, \nabla)$ be a Fedosov manifold admitting a metaplectic structure, $\mathcal{S} \to M$  
the corresponding symplectic spinor bundle and $d^{\nabla^S}$ the exterior covariant derivative. 
Then for each $\xi \in \Gamma(M,T^*M)$ and $\alpha \otimes \phi \in \Gamma(M, \mathcal{E}),$ 
the symbol $\sigma^{\xi}$ of $d^{\nabla^S}$  is given by 
$$\sigma^{\xi}(\alpha \otimes \phi):=  \xi \wedge \alpha \otimes \phi.$$

{\it Proof.} For $f \in \mathcal{C}^{\infty}(M),$ $\xi\in \Gamma(M,T^*M)$ and $\alpha \otimes s \in \Gamma(M,\mathcal{E}),$
let us compute $d^{\nabla^S}(f\alpha \otimes s) - f d^{\nabla^S}(\alpha \otimes s)=df \wedge \alpha \otimes s + fd^{\nabla^S}(\alpha \otimes s) 
- f d^{\nabla^S}(\alpha \otimes s) = df \wedge \alpha \otimes s.$
Using this computation, we get the statement of the lemma. $\Box$



Due to the previous lemma and the definition of the symplectic twistor operators,
we get that for each $i=0,\ldots, 2l$ and $\xi \in \Gamma(M,T^*M),$ the symbol $\sigma^{\xi}_i$ of 
the symplectic twistor operator $T_i$ is given by the formula
$$\sigma^{\xi}_i(\alpha \otimes s):= p^{i+1,m_{i+1}} (\xi \wedge \alpha \otimes s)$$
for each $\alpha \otimes s \in \Gamma(M, \mathcal{E}^i).$

In order to prove the ellipticity of the symplectic twistor complexes, i.e., the exactness of 
 the sequence of sheaves
 $$0  \longrightarrow
  \Gamma(M,\mathcal{E}^{0}) \overset{\sigma_0^{\xi}}{\longrightarrow} 
  \Gamma(M,\mathcal{E}^{1}) \overset{\sigma_1^{\xi}}{\longrightarrow}  
  \cdots  \overset{\sigma_{l-2}^{\xi}}{\longrightarrow}
  \Gamma(M,\mathcal{E}^{l-1})  \mbox{   and}  
$$ 
$$  \Gamma(M,\mathcal{E}^{l}) \overset{\sigma_l^{\xi}}{\longrightarrow} 
  \Gamma(M,\mathcal{E}^{l+1}) \overset{\sigma_{l+1}^{\xi}}{\longrightarrow}  
  \cdots  \overset{\sigma_{2l-1}^{\xi}}{\longrightarrow}
  \Gamma(M,\mathcal{E}^{2l}) \longrightarrow 0,  
$$ 
for any $\xi \in \Gamma(M,T^{*}M) \setminus \{0\},$  we need to compare 
the kernels and the images of the symbols maps. Therefore, we prove the following statement,
 in which the projections $p^{ii}$ are more specified. From now on, we   shall denote the projections $p^{ii}$ onto
 ${\bf E}^{i}$ by $p^i$ simply. (We will make no use  of the projections from ${\bf E}$ onto $\bigwedge^i\mathbb{V}^* \otimes {\bf S}$ or of their lifts to associated structures.)

{\bf Lemma 7:} For $i=0,\ldots,l-1,$ $\xi \in \mathbb{V}^*$ and $\alpha \otimes s \in {\bf E}^{i},$   we have 
\begin{eqnarray}
p^{i+1}(\xi \wedge \alpha \otimes s) = \xi \wedge \alpha \otimes s + \beta F^+ (\alpha \otimes \xi^{\sharp}. s) +\gamma (E^+ 
\iota_{\xi^{\sharp}}\alpha \otimes s) \label{p},
\end{eqnarray}
 where $\beta = \frac{2}{i-l} \mbox{ and  }\gamma = \frac{\imath}{i-l}.$

{\it Proof.} 
We split the proof into four parts.
\begin{itemize}

\item[1.] In this item, we prove that for a fixed $i \in \{0,\ldots, l\}$ and any $k=0,\ldots, i,$  
there exists $\alpha_k^i \in \mathbb{C}$ such that
$$p^{i}=\sum_{k=0}^{i}\alpha_k^i (F^+)^k (F^-)^k$$ with $\alpha_0^i=1$ for each $i=0,\ldots, l.$
Because for each $i=0,\ldots, l,$ the projections $p^{i}$ are $\tilde{G}$-equivariant morphisms, 
they can be expressed as (finite) linear combinations of the elements of the vector space $\mbox{End}_{\tilde{G}}({\bf E}).$ 
Due to the Lemma 3 item 4 (cf. also Kr\'ysl \cite{KryslJOLT2}), we know that the complex 
associative algebra $\mbox{End}_{\tilde{G}}({\bf E})$ is (finitely) generated by  $F^+$ and $F^-$  and the projections $p_{\pm}.$
It is easy to see that the projections $p_{\pm}$ can be omitted from any expression for $p^{i}$ and thus, 
each projection $p^{i}$  can be expressed just using $F^+$ and $F^-.$   Due to the defining relation
$H=2\{F^+, F^-\}$ and the relation (\ref{h}) on the values of $H$ on homogeneous elements, one can order the operators $F^+$ and $F^-$ in the expression for
$p^{i}$ and put the operators $F^+$ to the left-hand   and the operators  $F^-$ to the right-hand. In this way, we express $p^{i}$ as a linear combination of the expressions of type $(F^+)^a(F^-)^b$ for $a,b \in \mathbb{N}_0.$

Since the projection $p^{i}$ does not change the form degree of a symplectic spinor valued exterior form and $F^-$ and $F^+$ decreases and increases the form degree by one, respectively, the 
relation $a=b$ follows.  
Because the operator $F^-$  decreases the form degree by one, for $k>i,$ the summands $(F^+)^k(F^-)^k$  actually do not occur in the expression for the projection $p^{i}$ written above.
Thus, 
\begin{eqnarray}
p^{i}=\sum_{k=0}^i \alpha^i_k(F^+)^k (F^-)^k \label{pecko}
\end{eqnarray} 
for some $\alpha^i_k \in \mathbb{C},$ $k=0,\ldots, i.$

Now, we shall prove the equation $\alpha_0^i=1,$ $i=0,\ldots, l.$   By evaluating the left-hand side of (\ref{pecko}) on an 
element $\phi \in {\bf E}^{i}$ we get $\phi$, whereas at the right-hand side 
the only summand which remains is the zeroth one. (The other summands vanish because $F^-$ is $\tilde{G}$-equivariant, decreases the form degree by one 
and there is no summand in $\bigwedge^{i-1}\mathbb{V}^* \otimes {\bf S}$ isomorphic to ${\bf E}^{i}_+$ or 
to ${\bf E}^{i}_-.$ See the Remark item 3 bellow the Theorem 1.)

\item[2.] Now, suppose $\xi \in \mathbb{V}^*$ and $\alpha \otimes s \in {\bf E}^{i},$ $i=0,\ldots, l-1.$
 Due to the Theorem 2, we know that $\phi:=\xi \wedge \alpha \otimes s \in {\bf E}^{i+1,i-1} \oplus {\bf E}^{i+1,i} \oplus {\bf E}^{i+1,i+1}.$
Applying $p^{i+1}$ on the element $\phi,$ only the zeroth, first, and second summand in the expression 
$p^{i+1}\phi = \sum_{k=0}^{i+1}\alpha_k^{i+1} (F^+)^k (F^-)^k \phi.$ (For $ k > 2,$ the $k^{th}$ summand vanishes in the expression for $p^{i+1}\phi$  because $F^-$ is $\tilde{G}$-equivariant, 
decreases the form degree by one and there is no summand in $\bigwedge^{i-2} \mathbb{V}^* \otimes {\bf S}$ isomorphic to ${\bf E}^{i+1,i-1}_{\pm}$ 
or ${\bf E}^{i+1,i}_{\pm}$ or ${\bf E}^{i+1,i+1}_{\pm}.$ See the Remark bellow the Theorem 1.)

\item[3.] Due to the previous item, we already know that for the element $\phi = \xi \wedge \alpha \otimes s$ chosen above, we get $$p^{i+1}\phi =  \sum_{k=0}^2\alpha_k^{i+1} (F^{+})^k(F^-)^k\phi \label{formule}.$$ 
 Using the relations (\ref{formule}) and (\ref{ecko}), may write
\begin{eqnarray*}
&&p^{i+1}(\xi\wedge \alpha \otimes s)=\\
&=&\xi\wedge\alpha \otimes s + \alpha^{i+1}_1 F^+ F^-(\xi\wedge \alpha \otimes s)
+\alpha^{i+1}_2 (F^+)^2(F^-)^2(\xi\wedge \alpha \otimes s)\\
&=&\xi\wedge\alpha\otimes s +\alpha^{i+1}_1\frac{1}{2}F^+\omega^{ij}[(\iota_{e_i}\xi) \alpha \otimes e_j.s - \xi \wedge \iota_{e_i}\alpha \otimes e_j.s] -\\
&&\alpha^{i+1}_2 E^+ \frac{\imath}{32}\omega^{ij}\iota_{e_i}\iota_{e_j}(\xi\wedge \alpha \otimes s)\\
&=&\xi\wedge\alpha\otimes s - \alpha_1^{i+1}\frac{1}{2} F^+[\alpha \otimes \xi^{\sharp}.s + 2\xi \wedge F^-(\alpha\otimes s)]-\\
&&-\alpha^{i+1}_2E^+\frac{\imath}{32}\omega^{ij}\iota_{e_i}(\xi_j \alpha \otimes s 
-\xi\wedge\iota_{e_j}\alpha \otimes s)
\end{eqnarray*}
Because $\alpha \otimes s \in {\bf E}^{i},$ we get $F^-(\alpha\otimes s) = 0$ by Lemma 3 item 2.
Using this equation, we may write
\begin{eqnarray*}
p^{i+1}(\xi \wedge \alpha \otimes s)&=& \xi\wedge\alpha\otimes s -
\frac{\alpha_1^{i+1}}{2} F^+(\alpha \otimes \xi^{\sharp}.s) \\
&&- \frac{\imath \alpha_2^{i+1}}{32}E^+(2
   \xi^i \iota_{e_i}\alpha \otimes s 
 +\frac{2\alpha_2^{i+1}}{\imath}\xi \wedge E^-\alpha\otimes s).
\end{eqnarray*}
The last summand in this expression vanishes due to the Lemma 3 item 2 because  $E^-=-4F^-F^-$ and
$\alpha \otimes s \in {\bf E}^{i}.$ 

Summing-up, we have
$$p^{i+1}\phi = \xi\wedge\alpha\otimes s - \alpha_1^{i+1} \frac{1}{2} F^+(\alpha \otimes \xi^{\sharp}.s) -
\alpha_2^{i+1}\frac{\imath}{16} E^+\iota_{\xi^{\sharp}}\alpha \otimes s,$$
which is the formula of the form written in the statement of the lemma.

\item[4.] In this item, we shall determine the 
numbers $\beta, \gamma \in \mathbb{C}$. This can be done in at least two different ways.
The first one relies on the evaluation of $p^{i+1}$ (obtained in the preceding item) on an element $\phi \in {\bf E}^{i+1,i}.$ After a straightforward but tedious computation,
we get $\alpha_1^{i+1} = 4/(l-i)$ just using the formulas (\ref{h}), (\ref{ffo}) and (\ref{ffec}).
Evaluating $p^{i+1}$ on an element $\phi \in {\bf E}^{i+1,i-1},$ we get $\alpha_2^{i+1}=16/(l-i)$ using the same method. The second possibility is to use the identities and the idempotence $(p^{i+1})^2=p^{i+1}.$  
Thus, comparing the last written formula of the preceding item and the Eqn. (\ref{p}), we get $\beta =2/(i-l)$ and $\gamma = \imath/(i-l).$
\end{itemize} $\Box$

{\bf Remark:} 
For $i=l,\ldots, 2l,$ $\xi \in \mathbb{V}^{*} $ and $\alpha \otimes s \in {\bf E}^{i},$ the formula for $p^{i+1}$ reads simply
$$p^{i+1}(\xi \wedge \alpha \otimes s) = \xi \wedge \alpha \otimes s$$  because of the Theorem 2 and the Remark 3 bellow the Theorem 1. 







Now, using the Lemmas 3, 4, 6 a 7, we are able to prove the ellipticity of the symplectic twistor complexes.

{\bf Theorem 8:} Let $(M,\omega, \nabla)$ be a Fedosov manifold of Ricci type admitting a metaplectic structure. 
Then the truncated symplectic twistor complexes 
$$0  \longrightarrow
  \Gamma(M,\mathcal{E}^{0}) \overset{T_0}{\longrightarrow} 
  \Gamma(M,\mathcal{E}^{1}) \overset{T_{1}}{\longrightarrow}  
  \cdots  \overset{T_{l-2}}{\longrightarrow}  \Gamma(M,\mathcal{E}^{l-1}) \mbox{   and}  
$$ 
$$
  \Gamma(M,\mathcal{E}^{l}) \overset{T_l}{\longrightarrow} 
  \Gamma(M,\mathcal{E}^{l+1}) \overset{T_{l+1}}{\longrightarrow}  
  \cdots  \overset{T_{2l-1}}{\longrightarrow}
  \Gamma(M,\mathcal{E}^{2l}) \longrightarrow 0  
$$ 
 are elliptic. 

{\it Proof.} 

\begin{itemize}

\item[1.] First, we prove that the sequences mentioned in the formulation of the theorem are complexes. For
$i=0,\ldots,l-2, l, \ldots, 2l-1,$ $\psi \in \Gamma(M,\mathcal{E}^{i})$ and a differential 1-form $\xi \in \Gamma(M,T^*M),$   
we may write $0 = p^{i+1}(0) = p^{i+1}((\xi \wedge \xi) \wedge \psi)=
p^{i+1}(\xi \wedge \mbox{Id} (\xi \wedge \psi)) = p^{i+1}(\xi \wedge \sum_{j=0}^{m_{i+1}} p^{i+1,j}(\xi  \wedge \psi)).$ 
Due to the Theorem 2, we know that the last written expression equals
$p^{i+1}(\xi \wedge p^{i} (\xi \wedge \psi))= 
\sigma_{i+1}^{\xi}\sigma_i^{\xi}(\psi)$ and thus $\sigma_{i+1}^{\xi}\sigma_i^{\xi}=0.$ 
 
\item[2.] Second, we prove the relation $\mbox{Ker}(\sigma_{i}^{\xi}) \subseteq \overline{\mbox{Im}(\sigma_{i-1}^{\xi})}$\footnote{The closure is with respect to the Grothendieck tensor product topology.}
for each $0\neq \xi \in \mathfrak{X}(M)$ and $i=0,\ldots, l-2$ (recall $\sigma_{-1}^{\xi}=0$ due to our convention).
Suppose a homogeneous element $\alpha \otimes s \in \Gamma(M,\mathcal{E}^{i})$  is given such that $\sigma_{i}^{\xi}(\alpha \otimes s)=0.$ (At the end of this item, we will treat the general non-homogeneous case.)
Due to the paragraph bellow the Lemma 6 we know that $0=\sigma_{i}^{\xi}(\alpha\otimes s)= p^{i+1}(\xi \wedge \alpha \otimes s).$  
Thus, we shall find an element $\psi \in \Gamma(M, \mathcal{E}^{i-1})$ such that
$p^{i}(\xi \wedge \psi) = \alpha \otimes s.$ 

Using formula (\ref{p}) for the projection (Lemma 7), 
we may rewrite the equation $p^{i+1}(\xi \wedge \alpha \otimes s) = 0$ into 
\begin{eqnarray}
\xi \wedge \alpha \otimes s + \beta F^+ (\alpha \otimes \xi^{\sharp}.s) + \gamma E^+ \iota_{\xi^{\sharp}}\alpha \otimes s  = 0 \label{proj}.
\end{eqnarray}

  Applying the operator $E^-$ (in the form of the formula (\ref{ecko})) on the both sides of the previous equation 
and using the first commutation relation in the row (\ref{ef}) from Lemma 4,
we get
$$\frac{\imath}{2}\omega^{ij}\iota_{e_i}\iota_{e_j}(\xi \wedge\alpha) \otimes s 
+ \beta E^-F^+ ( \alpha \otimes \xi^{\sharp}. s) + $$
$$ +\gamma (E^+ E^- - 2 H)   \iota_{\xi^{\sharp}}\alpha \otimes s  = 0$$
 
Using the graded Leibniz property of $\iota_{\xi^{\sharp}},$ 
the relation (\ref{h}) for the values of $H$ on homogeneous elements and the second relation in the row (\ref{ef}) from Lemma 4, we obtain
$$\frac{\imath}{2}(-2\iota_{\xi^{\sharp}}  - 2\imath \xi \wedge E^-)(\alpha \otimes s) +  
\beta F^+E^-(\alpha \otimes \xi^{\sharp}.s) - \beta F^-(\alpha \otimes \xi^{\sharp}.s) + $$
$$ + \gamma E^+ E^- \iota_{\xi^{\sharp}} \alpha \otimes s 
+ \gamma(l-i+1)\iota_{\xi^{\sharp}}\alpha \otimes s =0.$$
Because $E^- = \frac{\imath}{2}\omega^{ij}\iota_{e_i}\iota_{e_j}$ (formula (\ref{ecko})), the operator $E^-$ 
commutes with the operator of the symplectic Clifford multiplication (by the vector field $\xi^{\sharp}$), and 
also with the contraction $\iota_{\xi^{\sharp}}.$ Using these two facts, we get
$$\frac{\imath}{2}(-2\iota_{\xi^{\sharp}} - 2\imath \xi \wedge  E^-)(\alpha \otimes s) + 
\beta F^+\xi^{\sharp}. E^-(\alpha \otimes s) - \beta F^-(\alpha\otimes \xi^{\sharp}.s) + $$
$$+ \gamma E^+ \iota_{\xi^{\sharp}} E^-\alpha \otimes s + \gamma(l-i+1)\iota_{\xi^{\sharp}}\alpha \otimes s = 0.$$

Because $F^-(\alpha \otimes s) = 0$ (Lemma 3 item 2) and thus, $E^-\alpha \otimes s (=  4 F^- F^- (\alpha \otimes s) ) = 0.$  Thus, we obtain the identity 
$$-\imath\iota_{\xi^{\sharp}}\alpha \otimes s  - \beta F^-(\alpha \otimes \xi^{\sharp}.s)
+ \gamma(l-i+1)\iota_{\xi^{\sharp}}\alpha \otimes s = 0.$$

Substituting the second relation in the row (\ref{fc}) into the previous equation and using the fact $F^-(\alpha \otimes s) = 0$ again, we get
$$-\imath\iota_{\xi^{\sharp}}\alpha \otimes s 
- \beta \xi^{\sharp}.F^- (\alpha \otimes s) - \beta\frac{\imath}{2} \iota_{\xi^{\sharp}}\alpha \otimes s+ $$
$$ +  \gamma(l-i+1)\iota_{\xi^{\sharp}}\alpha \otimes s  = 0.$$
Using the prescription for the numbers $\beta$ and $\gamma$ (Lemma 7) and the already twice used
relation $F^-(\alpha \otimes s) = 0,$ we get $(-\imath +\gamma(l-i+1) - \beta\frac{\imath}{2})\iota_{\xi^{\sharp}}
\alpha \otimes s = - 2 \imath \iota_{\xi^{\sharp}}\alpha \otimes s = 0$ from which
the equation 
\begin{eqnarray}
\iota_{\xi^{\sharp}}\alpha \otimes s  = 0 \label{kontr}
\end{eqnarray}
follows.  

Substituting this relation into the prescription for the projection $p^{i}$ (Eqn. (\ref{proj})),
we get for $i=0,\ldots, l-2$ the equation
 
\begin{eqnarray}
 0 &=& p^{i+1} (\xi \wedge \alpha \otimes s) = \xi \wedge \alpha \otimes s 
 + \beta F^+ (\alpha \otimes \xi^{\sharp}.s) \label{mezi}.
\end{eqnarray}

Applying the contraction operator $\iota_{\xi^{\sharp}}$  to the previous equation and using first formula in the row (\ref{fc}) from Lemma 4, we obtain

\begin{eqnarray*}
0 &=& -\xi \wedge \iota_{\xi^{\sharp}} \alpha \otimes s  - \beta F^+\iota_{\xi^{\sharp}}(\alpha \otimes \xi^{\sharp}.s) +\beta \frac{\imath}{2}\alpha \otimes \xi^{\sharp}.(\xi^{\sharp}.s).
\end{eqnarray*}

Using the fact that the contraction and symplectic Clifford multiplication commute, we get
\begin{eqnarray*}
0 &=& -\xi \wedge \iota_{\xi^{\sharp}} \alpha \otimes s  - \beta F^+ \xi^{\sharp}.(\iota_{\xi^{\sharp}}\alpha \otimes s) +\beta \frac{\imath}{2}\alpha \otimes \xi^{\sharp}.(\xi^{\sharp}.s).
\end{eqnarray*}

Substituting the Eqn. (\ref{kontr})  into the previous equation, we obtain
$$\alpha \otimes \xi^{\sharp}.(\xi^{\sharp}.s)  =  0.$$

Using the definition of $F^+,$ multiplying the equation (\ref{mezi})
by $\xi^{\sharp}$ and using the equation $\iota_{\xi^{\sharp}}\alpha \otimes s =0$ (Eqn. (\ref{kontr}) again), we get
\begin{eqnarray*}
0&=& \xi \wedge \alpha \otimes \xi^{\sharp}.s + \beta \frac{\imath}{2} \epsilon^i \wedge \alpha \otimes \xi^{\sharp}.e_i.\xi^{\sharp}.s,\\
0&=& \xi \wedge \alpha \otimes \xi^{\sharp}.s + \beta \frac{\imath}{2} \epsilon^i \wedge \alpha \otimes (e_i.\xi^{\sharp}.\xi^{\sharp}. - 
\imath \omega_0(\xi^{\sharp}, e_i)\xi^{\sharp}.)s.
\end{eqnarray*}

Substituting the identity $\alpha \otimes \xi^{\sharp}.\xi^{\sharp}.s=0$ into the previous equation, we obtain
$$0 =(1+\frac{1}{2}\beta)\xi \wedge \alpha \otimes \xi.^{\sharp}s.$$

If $i=0,\ldots, l-2,$ the coefficient $1+\beta/2 \neq 0,$ and thus by dividing,
we get  $\xi \wedge \alpha \otimes \xi^{\sharp}.s =0.$ 
Because the symplectic Clifford multiplication by a  non-zero vector is injective (see the subsection 2.2),
 we have $0=\xi\wedge\alpha \otimes s.$

All the computations in this item could have been done for a general $\phi \in \Gamma(M,\mathcal{E}^i)$ because of  linearity.  But for a better understanding, we were separating the form and the spinorial parts.
Thus, we have $\xi\wedge \phi=0.$
Let $\{\alpha^k\}_{k=1}^{r_i}$ be a local frame of exterior differential $i$ forms. We may write
$\phi = a_{kj}\alpha^k\otimes s^j$ for some smooth functions $a_{kj}$ on $M$ and local symplectic spinor fields $s^j,$ $k,j=1,\ldots, r_k.$
Because the operator $\xi \wedge$ acts only on the form-part of ${\bf E},$ we get
that $\xi\wedge a_{kj}\alpha^k = 0.$
Using the Cartan lemma on exterior forms, we get for each $j=1,\ldots, r_k,$  the existence of  exterior differential $i-1$ forms $\beta_j$ such that $a_{kj}\alpha^{k} = \xi \wedge \beta_j.$ 
Now, we can write
$a_{kj} \alpha^k \otimes s^j = p^{i}(a_{kj}\alpha^k \otimes s^j) = p^i(\xi \wedge \beta_j \otimes s^j) = 
\sigma_{i-1}^{\xi}(\beta_j \otimes s^j).$
Thus, we see that the element $\psi:=\beta_j \otimes s^j$ is the desired preimage of $\phi.$
 
\item[3.] Now,  we prove that $\mbox{Ker}(\sigma^{\xi}_i) \subseteq \overline{\mbox{Im} (\sigma_{i-1}^{\xi})}$
for $i = l+1, \ldots, 2l,$ $0\neq \xi \in \Gamma(M,T^*M).$  Let $\phi \in \mbox{Ker}(\sigma^{\xi}_i),$ then
$0= p^{i+1}(\xi\wedge \phi) =\xi\wedge \phi.$ As in the last paragraph of the previous item, we get
  the existence of a symplectic spinor valued exterior differential $(i-1)$ form $\psi$ such that $\phi = \xi \wedge \psi.$  
\end{itemize}
$\Box$


In the future, we would like to interpret the reduced cohomology groups,  
make further conclusions following from the proved ellipticity. Eventually, one can search for an application of symplectic twistor complexes in representation theory.

\end{document}